\documentclass[oneside,12pt]{amsart}

\usepackage{amscd,amsmath,latexsym,amssymb}
\usepackage[mathcal]{euscript}

\newtheorem{thm}{Theorem}[section]
\newtheorem{cor}[thm]{Corollary}
\newtheorem{lem}[thm]{Lemma}
\newtheorem{prop}[thm]{Proposition}

\theoremstyle{definition}
\newtheorem{defin}[thm]{Definition}
\theoremstyle{definition}

\newtheorem{example}[thm]{Example}
\newtheorem{remark}[thm]{Remark}
\newtheorem{conjecture}[thm]{Conjecture}

\newcommand{\Z}{{\mathbb Z}}

\begin{document}
\title{Remarks concerning Lubotzky's filtration}

\author[F.~R.~Cohen]{F.~R.~Cohen$^{*}$}
\address{Department of Mathematics,
University of Rochester, Rochester, NY 14627}
\email{cohf@math.rochester.edu}
\thanks{$^{*}$Partially supported by the NSF}

\author[Marston Conder]{Marston Conder}
\address{Department of Mathematics,
New Zealand Institute of Mathematics and its Applications}
\email{m.conder@auckland.ac.nz}

\author[J.~Lopez]{J.~Lopez}
\address{Department of Mathematics,
University of Rochester, Rochester, NY 14627}
\email{jlopez@math.rochester.edu}

\author[Stratos Prassidis]{Stratos Prassidis$^{**}$}
\address{Department of Mathematics
Canisius College, Buffalo, NY 14208, U.S.A.}
\email{prasside@canisius.edu}
\thanks{$^{**}$Partially supported by Canisius College Summer
Grant}

\begin{abstract}
A discrete group which admits a faithful, finite dimensional, linear
representation over a field $\mathbb F$ of characteristic zero is
called {\it linear}.  This note combines the natural structure of
semi-direct products with work of A.~Lubotzky \cite{lubotzky} on the
existence of linear representations to develop a technique to give
sufficient conditions to show that a semi-direct product is linear.

Let $G$ denote a discrete group which is a semi-direct product given
by a split extension
$$1 \to \pi \to G \to \Gamma \to 1.$$ This note defines an additional type
of structure for this semi-direct product called a {\it stable
extension } below. The main results are as follows:
\begin{enumerate}
\item If $\pi$ and $\Gamma$ are linear, and the extension is stable, then $G$ is
also linear. Restrictions concerning this extension are necessary to
guarantee that $G$ is linear as seen from properties of the
Formanek-Procesi ``poison group" \cite{formanek.procesi}.
\item If the action of $\Gamma$ on $\pi$ has a ``Galois-like" property that it
factors through the automorphisms of certain natural ``towers of
groups over $\pi$" ( to be defined below ), then the associated
extension is stable and thus $G$ is linear.
\item The condition of a stable extension also implies that $G$
admits filtration quotients which themselves give a natural
structure of Lie algebra and which also imply earlier results of
Kohno, and Falk-Randell \cite{kohno3,falk-rand} on the Lie algebra
attached to the descending central series associated to the
fundamental groups of complex hyperplane complements.
\end{enumerate}

The methods here suggest that a possible technique for obtaining new
linearity results may be to analyze automorphisms of towers of
groups.

\end{abstract}

\maketitle

\section{Introduction}

A. Lubotzky \cite{lubotzky}  or \cite{ddms}, pages 172-175, gave a
purely group theoretic criterion which is equivalent to the
existence a faithful finite dimensional representation over a field
$\mathbb F$ of characteristic zero for a discrete group $G$ (where
the image is not necessarily discrete). A group $G$ with this
property is called {\it linear}.

The purpose of this paper is to give an extension of Lubotzky's
criterion which can sometimes be applied to show that a semi-direct
product of linear groups is again linear. The main subject of this
article is a split extension of groups given by
\[
\begin{CD}
1 @>{}>> \pi @>{i}>>  G @>{p}>>  \Gamma @>{}>> 1
\end{CD}
\]
for which it is assumed that both $\pi$ and $\Gamma$ are linear. The
main purpose of this article is to define the notion of a {\it
stable extension} as given in Definition
\ref{defin:stable.Lubotzky.filtration} which implies that $G$ is
linear.

The approach weaves together semi-direct products regarded as
pull-backs of a universal semi-direct product known as the holomorph
together with certain choices of filtrations of both $\pi$ and
$\Gamma$. Roughly speaking, one of the main results here is that
representations of $\Gamma$ in the automorphism group of $\pi$ which
factor through the automorphism group of the filtration of $\pi$ as
given in Definition \ref{defin:automorphism.group.of.a.tower}
suffices to show that $G$ is linear via Lubotzky's criteria
\cite{lubotzky}.

Notice that it may be the case that both $\pi$ and $\Gamma$ admit
faithful finite dimensional representations, but that $G$ does not.
A basic example due to Formanek and Procesi \cite{formanek.procesi}
is a split extension
\[
\begin{CD}
1 @>{}>> F_3 @>{i}>> H @>{p}>>  F_2 @>{}>> 1
\end{CD}
\]
where $H = G $, $F_n$ is a free group on $n$ letters, and the group
$H$ admits the following presentation:
\begin{equation} \label{eq:poison pres}
H=\langle a_1,a_2,a_3,\phi_1,\phi_2 \mid \phi_i^{} a_j^{}
\phi_i^{-1}=a_j^{}, \phi_i^{} a_3^{} \phi_i^{-1}=a_3^{} a_i^{},
i,j=1,2\rangle.
\end{equation}
This example, the Formanek-Procesi ``poison group", is a subgroup of
$Aut(F_3)$, the automorphism group of $F_3$ and has the property
that the action of $F_2$ on the first homology group of $F_3$ is
non-trivial.

Contrasting examples with $\pi$ given by $F_n$ which do in fact
admit faithful finite dimensional representations from the methods
given here are explained next. A subgroup of $Aut(F_n)$ known as
McCool's group $M(n)$ is generated by automorphisms given by
conjugating a fixed basis element by another fixed basis element
\cite{mccool}. Furthermore, the kernel of the natural map $Aut(F_n)
\to GL(n,\mathbb Z)$, $IA_n$, contains $M(n)$.

Consider a split extension
\[
\begin{CD}
1 @>{}>> F_n @>{i}>> G @>{p}>>  \Gamma @>{}>> 1
\end{CD}
\]
where $\Gamma$ admits a faithful finite dimensional representation
and the action of $\Gamma$ on $F_n$ factors through $M(n)$. It is
shown below that $G$ is sometimes linear. Thus it is natural to ask
the following question which is also raised in \cite{ccp} with some
additional evidence here.

\begin{conjecture} \label{conjecture:faithful}
Consider a split exact sequence of groups
\[
\begin{CD}
1 @>>> F_n @>>> G @>>>  \Gamma @>>> 1
\end{CD}
\]
with $F_n$ a free group on $n$ letters and $\Gamma$ a group that
admits a finite dimensional faithful linear representation. If the
conjugation action of $\Gamma$ on $F_n$ is trivial on homology,
$H_1(F_n;\Z)$, and thus factors through $IA_n$, then $G$ is linear.

A weaker conjecture is to replace $IA_n$ by McCool's group $M(n)$:
that is, if the conjugation action of $\Gamma$ on $F_n$ factors
through $M(n)$, then $G$ is linear
\end{conjecture}
\

\begin{remark} \label{remark:faithful}
Observe that $M(2) = IA_2 $. Thus in case $n = 2$, this conjecture
follows directly from the observations in Corollary
\ref{cor:G.L.N.squared} below. In case $\Gamma$ is a subgroup of
$GL(n,\mathbb F)$, it follows from the computations below that $G$
is a subgroup of $GL(n + 4,\mathbb F)$ with details left as an
exercise.
\end{remark}


The authors would like to congratulate Tom Farrell and Lowell Jones
on this happy occasion of their 60-th birthday. The authors would
also like to thank the organizers for this stimulating and
interesting opportunity to participate in an excellent conference.

\section{Definitions and Statement of Results}\label{definitions}

Recall the following definition from \cite{ddms} (page 171) and
\cite{lubotzky}.

\begin{defin}\label{defin:filtration}
A {\it filtration of the group $\pi$} is a descending chain of
normal subgroups
$$\cdots \subseteq L_j(\pi)\subseteq \cdots
\subseteq L_1(\pi) \subseteq L_0(\pi)= \pi$$ for $j \geq 0$ such
that $\bigcap_ {j \geq 1} L_j(\pi) = \{1\}$.
\end{defin}

\begin{defin}\label{defin:p.conguence.system}
A {\it $p$-congruence system for the group $\pi$} is a filtration of
$\pi$ $$\cdots \subseteq L_j(\pi)\subseteq \cdots \subseteq L_1(\pi)
\subseteq L_0(\pi)= \pi$$ for $j \geq 0$ such that
\begin{enumerate}
  \item $\pi/ L_{1}(\pi)$ is finite, and
  \item $L_1(\pi)/ L_{1+j}(\pi)$ is a finite $p$-group for all $j\geq
  0$.
\end{enumerate}
\end{defin}

\begin{defin}\label{defin:original.Lubotzky.criterion}
A {\it bounded $p$-congruence system for the group $\pi$ } is a
$p$-congruence system for the group $\pi$ given by

  $$\cdots \subseteq L_j(\pi)\subseteq \cdots  \subseteq
L_1(\pi) \subseteq L_0(\pi)= \pi$$

\

\noindent such that $d(L_i(\pi)/ L_{j}(\pi)) \leq e$ for all $0 \leq
i < j$ (where the number $d(G)$ denotes the minimal number of
generators of the group $G$ (\cite{ddms}, page xvii)).  A bounded
p-congruence system is also called a {\bf Lubotzky filtration}
below.
\end{defin}

The following is a restatement here of a result due to A.~Lubotzky
\cite{lubotzky}.

\begin{thm}\label{prop:finite}
A  group $G$ admits a  bounded $p$-congruence system for some prime
$p$ if and only if $G$ admits a faithful finite dimensional
representation for some field of characteristic zero.
\end{thm}

Let $\mathrm{Aut}(\pi)$ denote the automorphism group of $\pi$.
Consider a discrete group $\pi$ together with the universal
semi-direct product $\mathrm{Hol}(\pi)$ ``the natural" split
extension of $\mathrm{Aut}(\pi)$ by $\pi$, $$1 \to\ \pi \to\
\mathrm{Hol}(\pi) \to\ \mathrm{Aut}(\pi) \to\ 1$$ The group
$\mathrm{Hol}(\pi)$, as a set, is the product $\mathrm{Aut}(\pi)
\times \pi$ with the product structure defined by the formula
$$(f,x)\cdot (g,y) = (f\cdot g, g^{-1}(x) \cdot y)$$ for $f,g$ in
$\mathrm{Aut}(\pi)$, and $x,y$ in $\pi$.

The next four formulas follow from the definition but are listed
here for convenience of the reader in the proofs below.

\begin{enumerate}
  \item $(f,x)^{-1} =  (f^{-1},f(x^{-1}))$,
  \item $(f,1)^{-1}\cdot (1,y)\cdot (f,1)= (1, f^{-1}(y))$,
  \item $(f,x)\cdot (g,y) \cdot (f,x)^{-1} = (f\cdot g\cdot f^{-1}, f(g^{-1}(x)\cdot y ) \cdot
  f(x^{-1}))$, and
  \item $[(f,x),(g,y)] = (f\cdot g\cdot f^{-1}\cdot g^{-1},g\{ f(g^{-1}(x)\cdot y ) \cdot f(x^{-1})\}\cdot g(y^{-1}))$.
\end{enumerate}

Consider a homomorphism $$\phi: \Gamma \to \mathrm{Aut}(\pi)$$
called the {\bf classifying map for the extension}. Pull back the
extension determined by $\mathrm{Hol}(\pi)$ to obtain the extension
$G$ together with a morphism of extensions (as developed in more
detail in \cite{voloshina}):
\[
\begin{CD}
1 @>>> \pi @>{i}>> G @>{p_{\phi}}>>  \Gamma @>>> 1\\
@. @VV{1}V          @VV{}V          @VV{\phi}V  @.\\
1 @>>> \pi @>{i}>> \mathrm{Hol}(\pi)   @>{p}>> \mathrm{Aut}(\pi)
@>>> 1
\end{CD}
\] Furthermore, every split extension of $\Gamma$ with kernel $\pi$ is
given by such a pull-back for some choice of homomorphism
$$\phi: \Gamma \to \mathrm{Aut}(\pi).$$
Thus if $(f,x), (g,y) \in \Gamma \times \pi$ then $(f,x)\cdot(g,y) =
(f\cdot g, \phi(g^{-1})(x)\cdot y)$. A notational convention used
throughout this article is that $g^{-1}(x)$ denotes
$\phi(g^{-1})(x)$.

\

The results here intertwine filtrations
for the groups $\pi$ and $\Gamma$ in the extension
\[
\begin{CD}
1 @>{}>> \pi @>{i}>>  G @>{p_{\phi}}>>  \Gamma @>{}>> 1
\end{CD}
\] by focusing on the classifying map for the extension
given by $\phi: \Gamma \to \mathrm{Aut}(\pi)$ rather than
considering the extension itself. Thus, the main focus here are
conditions concerning the homomorphism $\phi: \Gamma \to
\mathrm{Aut}(\pi)$ which imply that $G$ is linear.

\begin{defin}\label{defin:stable.Lubotzky.filtration}
Assume that
\[
\begin{CD}
1 @>{}>> \pi @>{i}>>  G @>{p_{\phi}}>>  \Gamma @>{}>> 1
\end{CD}
\]
is a split extension classified by the map $$\phi: \Gamma \to
\mathrm{Aut}(\pi)$$ together with filtrations
\begin{enumerate}
\item for the group $\pi$
$$\cdots \subseteq L_j(\pi)\subseteq \cdots  \subseteq L_1(\pi)
\subseteq L_0(\pi)= \pi$$ for $j \geq 0$  and
\item for $\Gamma$
$$\cdots \subseteq  F_j(\Gamma)\subseteq \cdots  \subseteq F_1(\Gamma)
\subseteq F_0(\Gamma)= \Gamma$$ for $j \geq 0$.
\end{enumerate}
{\it The extension ( together with the two
filtrations ) is said to be stable} if and only if for every $(g,y)$
in $ F_{r+s}(\Gamma) \times  L_{r+s}(\pi)$ and for every $(f,x)$ in
$ F_{r}(\Gamma) \times L_{r}(\pi)$ the following properties are
satisfied for $r, s \geq 0$:
\begin{enumerate}
  \item $f(y) \in  L_{r+s}(\pi)$ and
  \item $g(x)= \delta_x \cdot x$ for $\delta_x \in L_{r+s}(\pi)$.
\end{enumerate}
\end{defin}

\begin{remark}\label{remark:conditions}
These two conditions both of which must be satisfied in what is
given below fit naturally with extensions. They arise by considering
the natural ``twisting'' for the holomorph as well as for certain
fibre bundles.
\end{remark}


The definition of a stable extension is basically recording the
feature that certain extensions ``look like products" modulo certain
higher filtrations. One result is as follows.

\begin{thm} \label{thm:Lubotzky.filtrations}
Assume that the split extension
\[
\begin{CD}
1 @>{}>> \pi @>{i}>>  G @>{p_{\phi}}>>  \Gamma @>{}>> 1
\end{CD}
\] is classified by the map $\phi: \Gamma \to Aut(\pi)$ which
satisfies the conditions that
\begin{enumerate}
\item $\Gamma$ and $\pi$ admit bounded $p$-congruence systems for some prime $p$ as
given in Definition \ref{defin:original.Lubotzky.criterion}, and
\item the $p$-congruence systems for the groups $\pi$ and $\Gamma$ in
part $(1)$ are stable in the sense of Definition
\ref{defin:stable.Lubotzky.filtration}.
\end{enumerate}
Then $G$ is linear.
\end{thm}

Examples of Theorem \ref{thm:Lubotzky.filtrations} are given in
sections \ref{sec:An example} and \ref{sec:A second example}. These
examples arise by forming the split extension $$1 \to F_n \to G \to
\Gamma \to 1$$ where
\begin{enumerate}
  \item $\Gamma$ is a subgroup of $GL(2,\mathbb Z)$ ( and thus $\Gamma$ has
  a normal finite index subgroup which is free ),
  \item  $F_n$ is isomorphic to a principal congruence subgroup of
  level $p^r$  in $PSL(2,\mathbb Z)$, and
  \item $\Gamma$ acts by conjugation on $F_n$.
\end{enumerate} That these examples are linear follows from standard
elementary methods as well as the methods here. One related special
case is as follows.

\

\begin{example} \label{example:covering.extensions}
Consider the extension
\[
\begin{CD}
1 @>{}>>  F[a_1,a_2,\cdots,a_n,b] @>{}>>  G_n @>{}>> F[x,y]@>{}>> 1
\end{CD}
\] for which the action of $F[x,y]$ is given as follows.
        \begin{enumerate}
        \item \begin{enumerate}
        \item  $x(a_q) = a_{q+1}$ if $1 \leq q < n$
        with $x(a_n)= b\cdot a_1 \cdot b ^{-1}$ and
        \item  $x(b)= b$.
        \end{enumerate}

\item The action of $y$ is given by
        \begin{enumerate}
        \item $y(a_q) =  a_1 \cdot a_q \cdot a_1^{-1}$ and
          \item  $y(b)= a_1\cdot b \cdot  a_1^{-1}.$
        \end{enumerate}

\end{enumerate} Then $G_n$ is linear. As shown in section
\ref{sec:A second example}, these examples can be done easily by
using elementary, ``bare-hands" methods.
\end{example}

\

In the case of a split extension
\[
\begin{CD}
1 @>{}>> \pi @>{i}>>  G @>{p_{\phi}}>>  \Gamma @>{}>> 1
\end{CD}
\]
which is stable ( Definition \ref{defin:stable.Lubotzky.filtration}
), the group $G$ inherits a natural filtration which is defined next
with properties developed in section \ref{sec:Two.filtration}.

\begin{defin}\label{defin:two.filtration.extension}
A {\it filtration of the group $G$} is given by $$\mathbb F_j(G) =
F_j(\Gamma) \times L_j(\pi)$$ as a set with multiplication obtained
from restriction of the formula
$$(f,x)\cdot (g,y) = (f\cdot g, \phi(g)^{-1}(x) \cdot y)$$ for $f,g$
in $\Gamma$, and $x,y$ in $\pi$.
\end{defin}

\begin{remark}\label{remark:two.filtration.extension.remark}
To be precise, it must be checked that the stated multiplication in
Definition \ref{defin:two.filtration.extension} restricts to give
$\mathbb F_j(G)$ as a subgroup of $G$. This verification is carried
out in section \ref{sec:Two.filtration}.
\end{remark}

Let $H$ denote a discrete group. Recall that the commutator function
$$[-,-]:H \times H \to H$$ induces the structure of Lie algebra
on the associated graded for the descending central series
filtration of $H$. Kohno \cite{kohno3}, and Falk-Randell
\cite{falk-rand} obtained a structure theorem for these Lie algebras
restricted to certain semi-direct products of groups. A similar
theorem holds for the mod-$p$ descending central series filtration
\cite{cw}. However, there are other natural filtrations for which a
similar extension theorem holds which are addressed by using the
following definition.

\begin{defin}\label{defin:Lie.like.filtrations}
A filtration of the group $H$ given by $\{F_j(H)\}$ is said to be 
{\it Lie-like} provided the commutator function
$$[-,-]:H \times H \to H$$ restricts to $$[-,-]:F_p(H) \times F_q(H) \to F_{p+q}(H)$$ for all $ p, q \geq 0$.
\end{defin}

An analogue of this last property for split group extensions is
defined next.

\begin{defin}\label{defin:stably.lie.like}
Consider the split extension
\[
\begin{CD}
1 @>{}>> \pi @>{i}>>  G @>>>  \Gamma @>{}>> 1
\end{CD}
\]
Two filtrations $L_*(\pi)$ and $F_*(\Gamma)$ are said to be
{\it{stably Lie-like}} if
\begin{enumerate}
  \item $F_*(\Gamma)$ is Lie-like
  \item For $(f,x)\in F_r(\Gamma)\times L_r(\pi)$ and $(g,y)\in F_s(\Gamma)\times L_s(\pi)$, $f(x\cdot g(y))\in L_{r+s}(\pi)$
\end{enumerate}
\end{defin}

\begin{remark}
If $L_*(\pi)$ is a filtration as part of a stably Lie-like extension, then it is Lie-like.  For this, notice that $(1,x)\in F_r(\Gamma)\times L_r(\pi)$ and $(1,y)\in F_s(\Gamma)\times L_s(\pi)$ implies that $xy\in L_{r+s}(\pi)$.  Similary, $x^{-1}y^{-1}\in L_{r+s}(\pi)$.  Thus, the commutator $[x,y]=xyx^{-1}y^{-1}\in L_{r+s}(\pi)$.
\end{remark}

\begin{thm} \label{thm:Lie.algebra.extensions.for.stable.extensions}
Assume that the split extension
\[
\begin{CD}
1 @>{}>> \pi @>{i}>>  G @>{p_{\phi}}>>  \Gamma @>{}>> 1
\end{CD}
\] is classified by the map $\phi: \Gamma \to Aut(\pi)$ which
satisfies the following conditions:
\begin{enumerate}
\item The groups $\Gamma$ and $\pi$ admit filtrations (not
necessarily bounded $p$-congruence systems) $F_*(\Gamma)$ and
$L_*(\pi)$  as given in Definition \ref{defin:filtration}.
\item The filtrations for the groups $\pi$ and $\Gamma$ in
part $(1)$ are stable in the sense of Definition
\ref{defin:stable.Lubotzky.filtration}.
\item The filtrations $F_*(\Gamma)$ and $L_*(\pi)$ are both
stably Lie-like with associated graded Lie algebras denoted
$gr^F_*(\Gamma)$ and $gr^L_*(\pi)$.
\end{enumerate}

\noindent Then the filtration of $G$ given in Definition
\ref{defin:two.filtration.extension} is Lie-like.  Furthermore, there is a split, short exact sequence of Lie
algebras

$$0 \to gr^L_*(\pi) \to gr_*(G) \to gr^F_*(\Gamma) \to 0$$

\noindent where $gr_*(G)$ is the associated graded Lie algebra with Lie
bracket induced by the commutator pairing
$$[-,-]:G \times G \to G$$

\end{thm}

A systematic setting for stable extensions arises by considering
automorphisms of a tower of groups given by a bounded $p$-congruence
system for the group $\pi$. That method is recorded in the next
section.

\

\section{Automorphisms of Towers of Groups}\label{sec:Automorphisms of towers of groups}

The purpose of this section is (i) to define the automorphism group
of a tower of groups over a discrete group $\pi$ and (ii) to show
how the structure of the automorphism group of certain towers over
$\pi$ gives rise to linear groups. The automorphism group of a tower
of groups is defined next and is analogous to that of
\cite{Schneps}.

\begin{defin}\label{defin:automorphism.group.of.a.tower}
A {\it tower of groups over $\pi$} is

    \begin{enumerate}
  \item a family of groups $L_{n}(\pi)$ for $n$ in a pointed, totally ordered index set $I = S \cup \{\bullet\}$
  with unique least element  $\bullet$ and $L_{\bullet}(\pi) = \pi$,
  \item for every $i \geq j \in I$, there is a ( possibly empty ) family of
homomorphisms $\mathcal F(i,j)$ given by $\alpha(i,j):L_{i}(\pi)\to
L_{j}(\pi)$ with unique homomorphisms
$\alpha(i,\bullet):L_{i}(\pi)\to \pi$ such that
$$ \alpha(i,\bullet) =  \alpha(j,\bullet)\circ\alpha(i,j)$$ for all
$\alpha(i,j) \in \mathcal F(i,j)$ .

\end{enumerate}

\

The {\it automorphism group of this tower over $\pi$ denoted
$$Aut(L_*(\pi))$$ } is the subgroup of elements $(\phi_n) \in \prod_{n \in I}Aut(L_{n}(\pi))$
such that $$\phi_j\circ \alpha(i,j) = \alpha(i,j)\circ \phi_i, \;\; \text {for
all}\,\, \alpha(i,j) \in \mathcal F(i,j).$$

\end{defin}

\

A special case is given next.
\begin{defin}\label{defin:automorphism.of.an.inductive.tower}
An {\it inductive tower of groups over $\pi$} is a tower of groups
$\{L_{n}(\pi)\, | \, n \in I \, \cup \, \{\bullet\} \}$ over $\pi$ such that

\begin{enumerate}
\item the index set $I$ is given by the natural numbers $\mathbb N = I$ with $\bullet =
    0$,
\item each group $L_{n}(\pi)$ is a subgroup of $\pi$, and
\item for every $i\geq j$, there is exactly one $\alpha(i,j):L_{i}(\pi)\to
L_{j}(\pi)$ given by the natural inclusion.
\end{enumerate}

\end{defin}

Three remarks are given next.
\begin{remark}\label{remark:automorphism.D.C.S}

\begin{enumerate}
\item A filtration of $\pi$ given by $$\cdots \subseteq L_j(\pi)\subseteq \cdots  \subseteq L_1(\pi)
\subseteq L_0(\pi)= \pi$$ is an inductive tower over $\pi$. Thus, a
bounded $p$-congruence system is an inductive tower over $\pi$.

\item The automorphism group of an inductive tower of groups over $\pi$ is
the subgroup of elements in $Aut(\pi)$ which leave every
$L_{n}(\pi)$ invariant.

\item Restrict to the case where $L_{n}(\pi)$ is the
$(n+1)$-st stage of the descending central series of $\pi$,
$\Gamma^{n+1}(\pi)$. The natural inclusions
$$\cdots \subseteq \Gamma^{n+1}(\pi) \subseteq \cdots  \subseteq \Gamma^{2}(\pi)
\subseteq \Gamma^{1}(\pi) = \pi$$ specify an inductive tower over
$\pi$ for which each $\Gamma^{n+1}(\pi)$ is invariant. Thus, the
automorphism group of the inductive tower given by the descending
central series is equal to $Aut(\pi)$. Similar remarks apply to the
mod-$p$ descending central series of $\pi$.
\end{enumerate}

\end{remark}

\

The next Lemma is a remark which follows from the above definitions.

\begin{lem}\label{lem:automorphisms.of.towers}
Assume that $\{L_{n}(\pi)\,|\, n \in I\}$ is an inductive tower of
groups over $\pi$ so that the automorphism group of this tower,
$Aut(L_*(\pi))$, is a subgroup of $Aut(\pi)$. Given an automorphism
$\rho \in Aut(L_*(\pi)),$ there is the natural induced split
extension

\[
\begin{CD}
1 @>{}>> \pi @>{i}>>  G @>{p_{\rho}}>>  \Gamma @>{}>> 1
\end{CD}
\] classified by regarding $\rho \in Aut(\pi).$

\end{lem}

\

Automorphisms of certain towers then have implications for whether
extensions are linear.

\begin{thm}\label{thm:towers}
Consider the split extension

\[
\begin{CD}
1 @>{}>> \pi @>{}>>  G @>{}>>  \Gamma @>{}>> 1
\end{CD}
\]

\noindent and suppose that the following conditions are satisfied:

\begin{enumerate}
  \item The filtration $L_*(\pi)$ is a Lubotzky filtration for the group
  $\pi$.
  \item The extension is classfied by a map $\rho:\,\Gamma\longrightarrow Aut(L_*(\pi))$ where $Aut(L_*(\pi))\subseteq
Aut(\pi)$ is the automorphism group of the tower $L_*(\pi)$.
  \item There exists a Lubotzky filtration $F_*(\Gamma)$ for the group $\Gamma$ such
  that the filtrations $F_*(\Gamma)$ and $L_*(\pi)$ satisfy condition (2) in Definition \ref{defin:stable.Lubotzky.filtration}.
\end{enumerate}

\noindent Then $G$ is linear.

\end{thm}

\begin{proof}
It suffices to show that the extension is stable in the sense of
Definition \ref{defin:stable.Lubotzky.filtration}, since the result
will then follow from Theorem \ref{thm:Lubotzky.filtrations}

Suppose $(f,x)\in F_r(\Gamma)\times L_r(\pi)$ and $(g,y)\in
F_{r+s}(\Gamma)\times L_{r+s}(\pi)$ where $r,s\geq0$.  Since the
action of $\Gamma$ is tower-preserving and $y\in L_{r+s}(\pi)$, it
follows that $f(y)\in L_{r+s}(\pi)$ and the extension is stable.
\end{proof}

\begin{remark} \label{remark:non.split.extensions}
The constructions in this section give a method to extend the
techniques here to arbitrary group extensions without the assumption
that the extension is required to be split. This remark will be
addressed elsewhere.
\end{remark}

\section{Two Filtrations}\label{sec:Two.filtration}

The purpose of this section is to investigate split extensions
equipped with two filtrations as given in Definition
\ref{defin:stable.Lubotzky.filtration}. Suppose
\[
\begin{CD}
1 @>{}>> \pi @>{i}>>  G @>{p_{\phi}}>>  \Gamma @>{}>> 1
\end{CD}
\]
is a split extension classified by the map $\phi: \Gamma \to
\mathrm{Aut}(\pi)$ together with filtrations
\begin{enumerate}
\item $L_*(\pi)$ given by
$$ \cdots \subseteq L_j(\pi)\subseteq \cdots
\subseteq L_1(\pi) \subseteq L_0(\pi)= \pi$$ for $j \geq 0$ for the
group $\pi$ and
\item $F_*(\Gamma)$ given by
$$\cdots \subseteq  F_j(\Gamma)\subseteq \cdots  \subseteq
F_1(\Gamma) \subseteq F_0(\Gamma)= \Gamma$$ for $j \geq 0$ for the
group $\Gamma$.
\end{enumerate}

Assume that the extension (together with the two filtrations) is
stable as in Definition \ref{defin:stable.Lubotzky.filtration}. An
equivalent technical formulation for the definition of a stable
extension is stated next.  Although elementary, direct, and technical, 
the next lemma is checked here as the second condition listed is the 
one actually used in the proofs of the theorems below.

\begin{lem}\label{lem:reformulated.stable.extensions}
Assume that every $(g,y)$ in $ F_{r+s}(\Gamma) \times  L_{r+s}(\pi)$
and every $(f,x)$ in $ F_{r}(\Gamma) \times L_{r}(\pi)$. The
formulas given in Definition \ref{defin:stable.Lubotzky.filtration}
by

\begin{enumerate}
  \item $f(y) \in  L_{r+s}(\pi)$ and
  \item $g(x)= \delta_x \cdot x$ for $\delta_x \in L_{r+s}(\pi)$.
\end{enumerate} are equivalent to

\begin{enumerate}
  \item $f(y) \in  L_{r+s}(\pi)$ and
  \item $g^{-1}(x)\cdot x^{-1} \in L_{r+s}(\pi)$.
\end{enumerate}
\end{lem}

\begin{proof}
Assume that every $(g,y)$ in $ F_{r+s}(\Gamma) \times L_{r+s}(\pi)$
and every $(f,x)$ in $ F_{r}(\Gamma) \times L_{r}(\pi)$. It suffices
to check that $g(x)= \delta_x \cdot x$ for $\delta_x \in
L_{r+s}(\pi)$ if and only if $g^{-1}(x)\cdot x^{-1} \in
L_{r+s}(\pi)$.

\begin{enumerate}
  \item Assume that $g^{-1}(x)\cdot x^{-1} \in L_{r+s}(\pi)$, and so
$g^{-1}(x)\cdot x^{-1} = \epsilon_x \in L_{r+s}(\pi)$. Thus
$g(\epsilon_x) = {\delta_x}^{-1}  \in L_{r+s}(\pi)$ by setting f = g
and $y = \epsilon_x$. Thus $x\cdot g(x^{-1}) = {\delta_x}^{-1}.$
  \item Assume that  $g^{-1}(x)\cdot x^{-1} = {\epsilon_x}^{-1} \in
L_{r+s}(\pi)$. Apply $g$ to obtain $x \cdot g(x^{-1}) =
g(g^{-1}(x)\cdot x^{-1}) = g({\epsilon}^{-1}) \in L_{r+s}(\pi)$.
\end{enumerate}
\end{proof}

A filtration of $G$, $\mathbb F_{*}(G)$, was defined in Definition
\ref{defin:two.filtration.extension} without verifying that it is a
filtration, namely $\mathbb F_j(G)$ is naturally a subgroup of $G$.
This fact is recorded next.

\begin{lem}\label{lem:morphism}
Assume that
\[
\begin{CD}
1 @>{}>> \pi @>{i}>>  G @>{p_{\phi}}>>  \Gamma @>{}>> 1
\end{CD}
\]
is a split extension classified by the map $\phi: \Gamma \to
\mathrm{Aut}(\pi)$ and which is stable with respect to filtrations
$L_*(\pi) $ and $F_*(\Gamma)$. Then $\mathbb F_j(G)$ is a group
which is naturally a subgroup of $G$ and there is a morphism of
extensions
\[
\begin{CD}
L_j(\pi) @>{i}>> \mathbb F_j(G)@>{p_{\phi}}>>  F_j(\Gamma) \\
@V{inclusion}VV          @VV{}V          @VV{inclusion}V  \\
\pi @>{i}>> G @>{p_{\phi}}>>  \Gamma.
\end{CD}
\]

\end{lem}

\begin{proof}
It suffices to check that $\mathbb F_j(G)$ is closed with respect to
the product in $G$ given by $(f,x)\cdot (g,y) = (f\cdot g,
\phi(g)^{-1}(x) \cdot y)$ for $f,g$ in $\Gamma$, and $x,y$ in $\pi$
where, by convention, $$g(x)= \phi(g)(x).$$

Assume that $f,g$ are in $F_j(\Gamma)$, and that $x,y$ are in
$L_j(\pi)$. By the ``stability" condition in Definition
\ref{defin:stable.Lubotzky.filtration}, $\phi(g)^{-1}(x)$ is in
$L_j(\pi)$. Thus $\phi(g)^{-1}(x) \cdot y$ is in $L_j(\pi)$. The
lemma follows by inspection.
\end{proof}

Properties of the groups $\mathbb F_j(G)$ are recorded in the next
lemma.
\begin{lem} \label{lem:subgroup.for.two.filtrations}
Let
\[
\begin{CD}
1 @>{}>> \pi @>{i}>>  G @>{p_{\phi}}>>  \Gamma @>{}>> 1
\end{CD}
\] denote a split extension classified by a map $\phi: \Gamma \to Aut(\pi)$
and which is stable with respect to filtrations $L_*(\pi)$ and
$F_*(\Gamma)$. Let $\mathbb F_j(G)$ denote the groups defined
earlier.

Then there are morphisms of split extensions
\[
\begin{CD}
1 @>>> L_{r+s}(\pi) @>{i}>>  \mathbb F_{r+s}(G) @>{p_{\phi}}>>  F_{r+s}(\Gamma)  @>>> 1\\
@. @VV{}V          @VV{}V          @VV{1}V  @.\\
1 @>>> L_r(\pi) @>{i}>>  \mathbb F_r(G) @>{p_{\phi}}>>  F_r(\Gamma) @>>> 1\\
@. @VV{}V          @VV{}V          @VV{\phi}V  @.\\
1 @>>> \pi @>{i}>> \mathrm{Hol}(\pi) @>{p}>> Aut(\pi) @>>> 1
\end{CD}
\]
for every $s\geq 0$. Furthermore, $\mathbb F_{r+s}(G)$ is a normal
subgroup of $\mathbb F_r(G)$ and there is an extension
\[
\begin{CD}
1  @>{}>>  L_r(\pi)/L_{r+s}(\pi) @>{i}>>  \mathbb F_r(G)/\mathbb
F_{r+s}(G) @>{p_{\phi}}>>  F_r(\Gamma)/F_{r+s}(\Gamma)  @>{}>>  1.
\end{CD}
\]
Thus if $F_r(\Gamma)/F_{r+s}(\Gamma)$ is generated by $c$ elements
and $L_r(\pi)/L_{r+s}(\pi)$ is generated by $d$ elements, then
$\mathbb F_r(G)/\mathbb F_{r+s}(G)$ is generated by $c+d$ elements.

\end{lem}

\begin{proof}
In the proof below, recall the convention that $f(x)= \phi(f)(x)$
for $x \in \pi$, $f \in \Gamma$ and $\phi: \Gamma \to Aut(\pi)$.
Since the split extension
\[
\begin{CD}
1 @>{}>> \pi @>{i}>> G @>{p_{\phi}}>>  \Gamma @>{}>> 1
\end{CD}
\] is classified by a map $\phi: \Gamma \to Aut(\pi)$
which is stable with respect to filtrations $L_*(\pi)$ and
$F_*(\Gamma)$, there is a morphism of split extensions
\[
\begin{CD}
 1 @>>> L_{r+s}(\pi) @>{i}>>  \mathbb F_{r+s}(G) @>{p_{\phi}}>>  F_{r+s}(\Gamma)  @>>> 1\\
 @. @VV{}V          @VV{}V          @VV{1}V  @.\\
1 @>>> L_r(\pi) @>{i}>>  \mathbb F_r(G) @>{p_{\phi}}>>  F_r(\Gamma) @>>> 1\\
@. @VV{}V          @VV{}V          @VV{\phi}V  @.\\
1 @>>> \pi @>{i}>> \mathrm{Hol}(\pi) @>{p}>> Aut(\pi) @>>> 1
\end{CD}
\]
by Lemma \ref{lem:subgroup.for.two.filtrations}.

To check that $\mathbb F_{r+s}(G)$ is a normal subgroup of $\mathbb
F_r(G)$  for any $s\geq0$, let $(f,x)$ denote an element in
$F_r(\Gamma) \times L_r(\pi)$ and $(g,y)$ an element in
$F_{r+s}(\Gamma) \times L_{r+s}(\pi)$. Then
$$(f,x)\cdot (g,y) \cdot (f,x)^{-1} = (f\cdot g\cdot f^{-1},
f(g^{-1}(x)\cdot y ) \cdot f(x^{-1})).$$ Notice that

\begin{enumerate}
  \item $f\cdot g\cdot f^{-1}$ is in $F_{r+s}(\Gamma)$ since it's a normal
  subgroup of $F_{r}(\Gamma)$,
  \item $f(g^{-1}(x)\cdot y ) \cdot f(x^{-1}) = f(g^{-1}(x)) \cdot f(y)
\cdot f(x^{-1})$,
  \item $f(g^{-1}(x)) \cdot f(x^{-1})$ is in $L_{r+s}(\pi)$ by
  stability,
  \item $y$ is in $L_{r+s}(\pi)$ by assumption, thus $f(y)$ is in $L_{r+s}(\pi)$ by
  stability,
  \item $f(x) \cdot f(y)\cdot f(x^{-1})$ is in $L_{r+s}(\pi)$ by
  stability and
  \item $f(g^{-1}(x)) \cdot f(y)\cdot f(x^{-1}) = f(g^{-1}(x))\cdot f(x^{-1})\cdot f(x) \cdot f(y)\cdot
f(x^{-1})$ is in $L_{r+s}(\pi)$.
\item Thus $\mathbb F_{r+s}(G)$ is a normal subgroup of $\mathbb
F_{r}(G)$.
\end{enumerate}

Since $\mathbb F_{r+s}(G)$ is a normal subgroup of $\mathbb
F_{r}(G)$, there is a morphism of extensions
\[
\begin{CD}
1 @>>>  L_{r+s}(\pi) @>{i}>>  \mathbb F_{r+s}(G) @>{p_{\phi}}>>
F_{r+s}(\Gamma)  @>>> 1\\
 @. @VV{}V          @VV{}V          @VV{1}V @.\\
1 @>>> L_r(\pi) @>{i}>>  \mathbb F_r(G) @>{p_{\phi}}>>  F_r(\Gamma) @>>> 1\\
@.@VV{}V          @VV{}V          @VV{\phi}V  @.\\
1 @>>> L_r(\pi)/L_{r+s}(\pi) @>{i}>>  \mathbb F_r(G)/\mathbb
F_{r+s}(G) @>{}>> F_r(\Gamma)/F_{r+s}(\Gamma) @>>> 1.
\end{CD}
\]

Since $F_r(\Gamma )/F_{r+s}(\Gamma )$ is generated by $d$ elements,
the subgroup of $\mathbb F_r(G)/\mathbb F_{r+s}(G)$ generated by
lifts of these elements together with $c$ elements which generate
the kernel then generate the entire group. The lemma follows.
 \end{proof}

\section{Two Filtrations Continued: Proof of Theorem \ref{thm:Lubotzky.filtrations}}\label{sec:Two.Filtrations.continued.Lubotzky.filtration}

The purpose of this section is to describe properties of filtrations
arising in Section \ref{sec:Two.filtration} inspired by work of
A.~Lubotzky who gave a sufficient condition for the existence of a
finite dimensional faithful representation of a discrete group
\cite{lubotzky}. Lubotzky's filtration condition is changed below to
fit questions for an extension theorem.

Given filtrations for $\Gamma$ and $\pi$ which are stable
for the group extension
\[
\begin{CD}
1 @>{}>> \pi @>{i}>>  G @>{p_{\phi}}>>  \Gamma @>{}>> 1,
\end{CD}
\] there are naturally associated semi-direct products $\mathbb F_j(G)$ defined in
section \ref{sec:Two.filtration}.

Properties of the groups $\mathbb F_j(G)$ are recorded in the next
lemma.

\begin{lem}\label{lem:Luboztky.filtrations.from.two.filtrations}
Let
\[
\begin{CD}
1 @>{}>> \pi @>{i}>>  G @>{p_{\phi}}>>  \Gamma @>{}>> 1
\end{CD}
\] denote a split extension classified by a map $\phi: \Gamma \to Aut(\pi)$
and which is stable with respect to filtrations $L_*(\pi)$ and
$F_*(\Gamma)$ which are also assumed to be $p$-congruence systems.
Then $\mathbb F_*(G)$ is a $p$-congruence system for $G$.

\end{lem}

\begin{proof} To check that $\mathbb F_*(G)$ is $p$-congruence system for
$G$, recall that it suffices to check (by Definition
\ref{defin:p.conguence.system}) that
$$\bigcap_{j \geq 1}\mathbb F_j(G) = \{1\}$$
and $\mathbb F_*(G)$ is a descending chain of normal subgroups
$$\cdots \subseteq \mathbb F_j(G)\subseteq \cdots
\subseteq \mathbb F_1(G) \subseteq \mathbb F_0(G) = G$$ for $j \geq
0$ such that
\begin{enumerate}
  \item $G/ \mathbb F_1(G)$ is finite and
  \item $\mathbb F_1(G)/ \mathbb F_{1+j}(G)$ is a finite $p$-group for all $j\geq
  0$.
\end{enumerate}

That $\mathbb F_{r+s}(G)$ is a normal subgroup of $\mathbb F_r(G)$
is checked in Lemma \ref{lem:subgroup.for.two.filtrations}. Notice
that by the proof of Lemma \ref{lem:subgroup.for.two.filtrations},
$$\bigcap_{j \geq 1}\mathbb F_j(G) = \bigcap_{j \geq 1} (F_j(\Gamma) \times
L_j(\pi)) = \{1\}.$$ Furthermore by
\ref{lem:subgroup.for.two.filtrations}, $\mathbb F_*(G)$ is a
decreasing filtration of $G$ with the property that there is an
extension
\[
\begin{CD}
1  @>{}>>  L_r(\pi)/L_{r+s}(\pi) @>{i}>>  \mathbb F_r(G)/\mathbb
F_{r+s}(G) @>{p_{\phi}}>>  F_r(\Gamma)/F_{r+s}(\Gamma)  @>{}>>  1.
\end{CD}
\]
Thus
\begin{enumerate}
  \item if $\Gamma/F_{j}(\Gamma)$ as well as $ \pi/L_{j}(\pi)$ are finite,
    then so is $G/\mathbb F_{j}(G)$ and
  \item if $F_r(\Gamma)/F_{r+s}(\Gamma)$ as
    well as $L_r(\pi)/L_{r+s}(\pi)$ are finite $p$-groups, then so
 is $\mathbb F_r(G)/\mathbb F_{r+s}(G)$.
\end{enumerate}

Thus $\mathbb F_*(G)$ is $p$-congruence system for $G$ and the lemma
follows.

\end{proof}

\begin{lem}\label{lem:Luboztky.criterion.from.two.filtrations}
Let
\[
\begin{CD}
1 @>{}>> \pi @>{i}>>  G @>{p_{\phi}}>>  \Gamma @>{}>> 1
\end{CD}
\] denote a split extension classified by a map $\phi: \Gamma \to Aut(\pi)$
and which is stable with respect to filtrations $L_*(\pi)$ and
$F_*(\Gamma)$ which are also assumed to be bounded $p$-congruence
systems. Then $\mathbb F_*(G)$ is a Lubotzky filtration, a bounded
$p$-congruence system.

\end{lem}

\begin{proof} By Lemma \ref{lem:Luboztky.filtrations.from.two.filtrations},
$\mathbb F_*(G)$ is $p$-congruence system for $G$. Furthermore by
\ref{lem:subgroup.for.two.filtrations}, $\mathbb F_*(G)$ is a
decreasing filtration of $G$ with the property that there is an
extension

\[
\begin{CD}
1  @>{}>>  L_r(\pi)/L_{r+s}(\pi) @>{i}>>  \mathbb F_r(G)/\mathbb
F_{r+s}(G) @>{p_{\phi}}>>  F_r(\Gamma)/F_{r+s}(\Gamma)  @>{}>>  1.
\end{CD}
\]

Thus if $F_r(\Gamma)/F_{r+s}(\Gamma)$ is generated by $c$ elements
and $L_r(\pi)/L_{r+s}(\pi)$ is generated by $d$ elements, then
$\mathbb F_r(G)/\mathbb F_{r+s}(G)$ is generated by $c+d$ elements.
By Definition \ref{defin:original.Lubotzky.criterion}, $\mathbb
F_*(G)$ is a a Lubotzky filtration, a bounded $p$-congruence system
for $G$.

\end{proof}

One consequence of Theorem \ref{prop:finite} as well as Lemma
\ref{lem:Luboztky.filtrations.from.two.filtrations} is Theorem
\ref{thm:Lubotzky.filtrations}.

\section{Proof of Theorem \ref{thm:Lie.algebra.extensions.for.stable.extensions}}\label{sec:Lie.algebra.extnsions.for.stable.extensions}

Consider the two filtrations $L_*(\pi)$ and $F_*(\Gamma)$ associated
to the stable extension

\[
\begin{CD}
1 @>{}>> \pi @>{i}>>  G @>{p_{\phi}}>>  \Gamma @>{}>> 1
\end{CD}
\] as stated in Definition \ref{defin:stable.Lubotzky.filtration} and
developed in Section \ref{sec:Two.filtration}.

A filtration of $G$ regarded as a set was defined by $$\mathbb
F_j(G) =  F_j(\Gamma) \times L_j(\pi)$$ in Definition
\ref{defin:two.filtration.extension}. Some properties of $\mathbb
F_j(G)$ were proven in Lemmas \ref{lem:morphism} and
\ref{lem:subgroup.for.two.filtrations} as follows.

\begin{enumerate}
  \item The subset $\mathbb F_j(G)$ is naturally a subgroup of $G$.
  \item There is a morphism of split group extensions
  \[
\begin{CD}
1 @>>> L_{r+s}(\pi) @>{i}>>  \mathbb F_{r+s}(G) @>{p_{\phi}}>>  F_{r+s}(\Gamma)  @>>> 1\\
@. @VV{}V          @VV{}V          @VV{1}V  @.\\
1 @>>> L_r(\pi) @>{i}>>  \mathbb F_r(G) @>{p_{\phi}}>>  F_r(\Gamma) @>>> 1.\\
\end{CD}
\]
\item There is a split extension
\[
\begin{CD}
1  @>{}>>  L_r(\pi)/L_{r+s}(\pi) @>{i}>>  \mathbb F_r(G)/\mathbb
F_{r+s}(G) @>{p_{\phi}}>>  F_r(\Gamma)/F_{r+s}(\Gamma)  @>{}>>  1.
\end{CD}
\]
\end{enumerate}

Consider the filtration quotients

$$gr_r^F(\Gamma) = F_r(\Gamma)/F_{r+1}(\Gamma), \quad
gr_r^L(\pi) = L_r(\pi)/L_{r+1}(\pi), \quad \text{and}\quad
gr_r(G) = \mathbb F_r(G)/\mathbb F_{r+1}(G).$$

Then there is a split short exact sequence of groups
\[
\begin{CD}
\{0\}  @>{}>> gr_r^L(\pi) @>{}>> gr_r(G)  @>{ }>> gr_r^F(\Gamma)
@>{}>> \{0\}
\end{CD}
\] by Lemmas \ref{lem:morphism} and \ref{lem:subgroup.for.two.filtrations}.

That $\mathbb{F}_*(G)$ is Lie-like is checked next.  Suppose
$(f,x)\in\mathbb{F}_s(G)$ and $(g,y)\in\mathbb{F}_r(G)$. It will be
checked that $g[f(g^{-1}(x)\cdot y)\cdot f(x^{-1})]\cdot
g(y^{-1})\in L_{r+s}(\pi)$ whenever the following conditions are
satisfied:

\begin{enumerate}
  \item The extenstion is stable.
  \item The filtration on $\Gamma$ is Lie-like.
  \item $f(x\cdot g(y))\in L_{r+s}(\pi)$.
\end{enumerate}

\vskip .1in

Since the filtration $F_*(\Gamma)$ is Lie-like, there exists $h\in
F_{r+s}(\Gamma)$ with $gfg^{-1}=fh$. Since the extension is stable,
there exists $\delta_x\in L_{r+s}(\pi)$ such that
$h(x)=\delta_x\cdot x$.  This implies the following:

$$\begin{array}{rcl}
g[f(g^{-1}(x)\cdot y)\cdot f(x^{-1})]\cdot g(y^{-1}) & = & gfg^{-1}(x)\cdot gf(y)\cdot gf(x^{-1})\cdot g(y^{-1}) \\[1em]
                                                     & = & fh(x)\cdot gf(y)\cdot gf(x^{-1})\cdot g(y^{-1}) \\[1em]
                                                     & = & f(\delta_x)\cdot f(x)\cdot gf(y)\cdot gf(x^{-1})\cdot g(y^{-1}) \\[1em]
\end{array}$$

\vskip .1in

\noindent Notice that $f(\delta_x)\in L_{r+s}(\pi)$ by stability. So
it suffices to show $f(x)\cdot gf(y)\cdot gf(x^{-1})\cdot
g(y^{-1})\in L_{r+s}(\pi)$. Since $F_*(\Gamma)$ is Lie-like, there
exists $k\in F_{r+s}(\Gamma)$ with $gf=fgk$.  Since the extension is
stable and the filtration of $\pi$ is given by normal subgroups,
there is $\delta_y\in L_{r+s}(\pi)$ such that $k(y)=y\cdot\delta_y$.
This implies the following:

$$\begin{array}{rcl}
f(x)\cdot gf(y)\cdot gf(x^{-1})\cdot g(y^{-1}) & = & f(x)\cdot fgk(y)\cdot gf(x^{-1})\cdot g(y^{-1}) \\[1em]
                                               & = & f(x)\cdot fg(y\cdot\delta_y)\cdot g(f(x^{-1})\cdot y^{-1}) \\[1em]
                                               & = & f(x)\cdot fg(y)\cdot fg(\delta_y)\cdot g(f(x^{-1})\cdot y^{-1}) \\[1em]
                                               & = & f(x\cdot g(y))\cdot fg(\delta_y)\cdot g(f(x^{-1})\cdot y^{-1}) \\[1em]
\end{array}$$

\vskip .1in

Now $fg(\delta_y)\in L_{r+s}(\pi)$ by the stability condition. The
additional condition (3) above gives that $f(x\cdot g(y))\in L_{r+s}(\pi)$
and $g(f(x^{-1})\cdot y^{-1})\in L_{r+s}(\pi)$.

\

To finish the proof, notice that Theorem
\ref{thm:Lie.algebra.extensions.for.stable.extensions} follows at
once from the property that these maps induce morphisms of Lie
algebras, a property which is checked next.

\

First observe that if $x \in L_r(\pi)$ and  $y \in L_s(\pi)$, then
$[x,y] \in L_{r+s}(\pi)$ by the assumption that the filtration
$L_*(\pi)$ is Lie-like. Secondly, since the filtration of $G$ is
Lie-like, there is a commutative diagram

\[
\begin{CD}
L_r(\pi) \times L_s(\pi) @>{[-,-]}>> L_{r+s}(\pi)\\
@V{i \times i}VV              @VV{i}V  \\
\mathbb{F}_r(G) \times \mathbb{F}_s(G) @>{[-,-]}>> \mathbb{F}_{r+s}(G).
\end{CD}
\]

Thus the map $i:  \pi \to G $ passes to quotients on the level of
associated graded modules and preserves the structure of the
underlying Lie algebras.  Thus the map $$p_{\phi}:  G
\to \Gamma$$ preserves the structure of Lie algebras. The Theorem
follows.

\section{An Example}\label{sec:An example}

The purpose of this section is to give examples of Theorem
\ref{thm:Lubotzky.filtrations} and Theorem \ref{thm:towers}. This
example has the serious drawback that the extension can be shown to
be linear by a ``bare-hands", more general, classical argument which
is reviewed in Section \ref{sec:A second example}.

\

These examples arise by forming the split extension $$1 \to F_n \to
G \to \Gamma \to 1$$ where
\begin{enumerate}
  \item $\Gamma$ is a subgroup of $GL(2,\mathbb Z)$ ( and thus $\Gamma$ has
  a normal finite index subgroup which is free ),
  \item  $F_n$ is isomorphic to a principal congruence subgroup of
  level $p^r$  in $PSL(2,\mathbb Z)$, and
  \item $\Gamma$ acts by conjugation on $F_n$.
\end{enumerate}

\


\

Let $P\Gamma(2,p^r)$ denote the kernel of the ``mod-$p^r$ reduction
map" $$ \rho_{p^r}: PSL(2, \mathbb Z) \longrightarrow PSL(2,
\mathbb Z/p^r \mathbb Z).$$ Natural automorphisms of
$P\Gamma(2,p^r)$ as well as the tower
$$\cdots \subseteq  P\Gamma(2,p^{r+1}) \subseteq P\Gamma(2,p^r) \subseteq \cdots  \subseteq P\Gamma(2,p) = \pi$$ are
given by conjugation by an element in $GL(2,\mathbb Z)$.

\

Furthermore if $p$ is a prime, the groups $P\Gamma(2,p)$ are free on
$1+p(p^2-1)/12$ generators if $p$ is an odd prime \cite{frasch} or
$2$ letters if $p=2$ \cite{f.klein}.  Let $\Gamma(2,p^r)$ denote the kernel of the natural reduction map $GL(2,\mathbb{Z})\to GL(2,\mathbb{Z}/p^r\mathbb{Z})$.  Below it is shown that 

$$\cdots P\Gamma(2,p^{r+1})\subseteq P\Gamma(2,p^r)\subseteq\cdots\subseteq P\Gamma(2,p^2)\subseteq P\Gamma(2,p)$$

\noindent gives a Lubotzky filtration for $P\Gamma(2,p)$.  The reader can check that similar arguments show that 

$$\cdots\subseteq\Gamma(2,p^{r+1})\subseteq\Gamma(2,p^r)\subseteq\cdots\subseteq\Gamma(2,p^2)\subseteq\Gamma(2,p)$$

\noindent gives a Lubotzky filtration for $\Gamma(2,p)$.  This information is recorded
next while a more standard development is given in Section
\ref{sec:A second example}.

\begin{lem}
The filtration 

$$\cdots\subseteq\Gamma(2,p^{r+1})\subseteq\Gamma(2,p^r)\subseteq\cdots\subseteq\Gamma(2,p^2)\subseteq\Gamma(2,p)$$ 

\noindent of $\Gamma(2,p)$ is a Lubotzky filtration.
\end{lem}


\begin{prop}\label{thm:an.example.of.the.method}
The extension

\[
\begin{CD}
1 @>{}>> P\Gamma(2,p^r) @>{}>>  G @>{}>>  \Gamma(2,p^s) @>{}>> 1
\end{CD}
\]

\noindent is linear where $\Gamma(2,p^s)$ acts on $P\Gamma(2,p^r)$
by conjugation and $r,s\geq1$.

\end{prop}

\begin{proof}
Let $f\in\Gamma(2,p^s)$ and $x\in P\Gamma(2,p^{r+q})$ where
$q\geq0$, so that $x$ projects to the identity in $PSL(2,\mathbb Z/
p^{r+q}\mathbb Z)$. Since $fxf^{-1}\in P\Gamma(2,p^{r+q})$, the
conjugation action is tower-preserving. Thus this filtration, along
with the filtration of $P\Gamma(2,p^r)$ is stable in the sense of
Definition \ref{defin:stable.Lubotzky.filtration}. Theorem
\ref{thm:Lubotzky.filtrations} then implies that $G$ is linear.

\end{proof}

\

Additional properties, some classical, some possibly not, are
recorded next. Notice that an automorphism of the tower of groups
$$\cdots \subseteq  P\Gamma(2,p^{r+1}) \subseteq P\Gamma(2,p^r) \subseteq \cdots  \subseteq P\Gamma(2,p) = \pi$$
induces an automorphism of the Lie algebra $$gr_*(P\Gamma(2,p))
=\oplus_{s \geq 1}P\Gamma(2,p^{s})/P\Gamma(2,p^{s+1}).$$ Thus it is
natural to identify the structure of this Lie algebra.

\

That structure is given next where related, and standard properties
of these principal congruence subgroups are recorded for
convenience. Recall that $P\Gamma(2,p^{s+1})$ is a normal subgroup
of $P\Gamma(2,p^{s})$. Define $$gr_s( P\Gamma(2,p))=
P\Gamma(2,p^{s})/P\Gamma(2,p^{s+1})$$ the associated graded.

\

The commutator map $$[-,-]:PSL(2, \mathbb Z) \times PSL(2, \mathbb
Z)  \to PSL(2, \mathbb Z) $$ restricts to
$$[-,-]:P\Gamma(2,p^r) \times P\Gamma(2,p^s) \to
P\Gamma(2,p^{r+s}),$$ and induces the structure of Lie algebra on
the associated graded $$gr_*(P\Gamma(2,p))) = \oplus_{s \geq 1}
gr_s(P\Gamma(2,p))$$ with $$[-,-]: gr_s( P\Gamma(2,p)) \otimes gr_t(
P\Gamma(2,p)) \to gr_{s+t}( P\Gamma(2,p)).$$

\

Furthermore, the $p$-th power map
$$\psi^p: P\Gamma(2,p^s) \to P\Gamma(2,p^{s+1})$$
induces a ( possibly non-linear ) map $$\psi^p: gr_s( P\Gamma(2,p))
\to gr_{s+1}( P\Gamma(2,p)).$$ Together with the previous structure
of Lie algebra for $gr_*(P\Gamma(2,p))$, this gives the structure of a
restricted Lie algebra over the field with $p$ elements $\mathbb
F_p$. Classical, well-known properties of the fitration quotients
$P\Gamma(2,p^r)/P\Gamma(2,p^{r+1})$ are recorded in the next
theorem.
\begin{thm}\label{thm:filtration.quotients.PGamma(2.p)}

If $p$ is an odd prime, there are isomorphisms
$$\theta_q: \oplus_3 \mathbb Z / p\mathbb Z \to gr_q(
P\Gamma(2,p))$$ with a choice of basis given by

 $$ A_q =
  \left(  \begin{array}{ccc}
            1 & p^q \\
            0 & 1 \\
          \end{array}
        \right), \quad
 B_q =
  \left(  \begin{array}{ccc}
            1 & 0 \\
            p^q & 1 \\
          \end{array}
        \right), \quad
 C_q =
  \left(  \begin{array}{ccc}
            1+p^q & p^q \\
            -p^q & 1-p^q \\
          \end{array}
        \right).
$$

Furthermore,

$$B_q\cdot C_q \cdot A_q^{-1} = D_q$$
( where the next matrix is not the reduction of a matrix in
$PSL(2,\mathbb Z)$ but represents a nontrivial coset in
$gr_q(P\Gamma(2,p))$ )

 $$D_q =
  \left(  \begin{array}{ccc}
            1+p^q & 0 \\
            0 & 1-p^q \\
          \end{array}
        \right).
$$

If $p =2$, there are isomorphisms
$$\theta_q: \oplus_2 \mathbb Z / 2\mathbb Z \to gr_q(
P\Gamma(2,2))$$ with a choice of basis given by

 $$ A_q =
  \left(  \begin{array}{ccc}
            1 & 2^q \\
            0 & 1 \\
          \end{array}
        \right), \quad
  B_q =
  \left(  \begin{array}{ccc}
            1 & 0 \\
            2^q & 1\\
          \end{array}
        \right).
$$
Furthermore $$[A_q,B_q] = 1.$$
\end{thm}

The additive structure given above is given in a global way in terns
of restricted Lie algebras. That structure is listed next.

\begin{thm}\label{thm:restricted.lie.algebras.for.PGamma(2.p)}

If $p=2$, then the restricted Lie algebra
$gr_*(P\Gamma(2,2))$ is generated by $A_1$ and $B_1$
(as a restricted Lie algebra). Furthermore, $gr_*(P\Gamma(2,2))$ is
the abelian, free restricted Lie algebra (over $\mathbb F_2$)
generated by $A_1$ and $B_1$ where, redundantly, the following
relations are satisfied:

\begin{enumerate}
  \item $[A_q,B_s] = 1$ for all $q$ and $s$,
\item $\psi^2(A_q) = A_{q+1}$ and
\item $\psi^2(B_q) = B_{q+1}$.
\end{enumerate}

If $p$ is an odd prime, then the restricted Lie algebra
$gr_*(P\Gamma(2,p))$ is generated by $A_1$, $B_1$ and $D_1$.
Furthermore, $gr_*(P\Gamma(2,p))$ is the free restricted Lie algebra
(over $\mathbb F_p$) generated by $A_1$, $B_1$ and $D_1$ subject to
the following relations.

\begin{enumerate}
  \item $[A_q,B_s] = D_{q+s}$ for all $q$ and $s$,
  \item $[A_q,D_s] = A_{q+s}^{-2}$ for all $q$ and $s$,
  \item $[B_q,D_s] = B_{q+s}^2$ for all $q$ and $s$,
\item $\psi^p(A_q) = A_{q+1}$,
\item $\psi^p(B_q) = B_{q+1}$, and
\item $\psi^p(D_q) = D_{q+1}$.
\end{enumerate}

\end{thm}

\

Theorem \ref{thm:filtration.quotients.PGamma(2.p)} is classical and
can be found in \cite{frasch}. The proof of Theorem
\ref{thm:restricted.lie.algebras.for.PGamma(2.p)} is a computation
based on the next classical lemma.

\begin{prop} \label{prop:extensions.for.congruence.subgroups.in.PSL(2.Z)}

The quotient $P\Gamma(2, p^r)/P\Gamma(2, p^{r+1})$ is isomorphic to
the kernel of the natural reduction map $$\gamma_{p^r}:PSL(2,
\mathbb Z/p^{r+1} \mathbb Z) \longrightarrow PSL(2, \mathbb Z/p^r
\mathbb Z)$$ and so there are isomorphisms

\[
P\Gamma(2, p^r)/P\Gamma(2, p^{r+1}) \cong
\begin{cases}
\oplus_2 \mathbb Z / 2\mathbb Z & \text{if $p = 2$, and}\\
\oplus_3 \mathbb Z / p\mathbb Z & \text{if $p$ is an odd prime.}
\end{cases}
\]

\end{prop}

\begin{proof}
The proof follows directly from the commutative diagram where the
rows and columns are all group extensions:

\[
\begin{CD}
P\Gamma(2, p^{r+1}) @>{1}>>  P\Gamma(2, p^{r+1}) @>{}>> 1\\
 @VV{}V          @VV{}V          @VV{1}V \\
P\Gamma(2, p^r)  @>{}>> PSL(2, \mathbb Z) @>{\rho_{p^r}}>> PSL(2, \mathbb Z/p^r \mathbb Z)  \\
@VV{}V          @VV{\rho_{p^{r+1}}}V          @VV{1}V  \\
Ker(2,p^r) @>{}>> PSL(2, \mathbb Z/p^{r+1} \mathbb Z)
@>{\gamma_{p^r}}>> PSL(2, \mathbb Z/p^r \mathbb Z)
\end{CD}
\]

\end{proof}

More applications to other groups $SL(n,A)$ and to their cohomology
will appear in the thesis of J.~Lopez \cite{j.lopez}.

\section{A Second Example}\label{sec:A second example}

The purpose of this section is to review classical properties of the
natural extension of $PSL(n,A)$ by $SL(n,A)$ with conjugation action
where $A$ is a commutative ring.  First, Let $Z(G)$ denote the
center of the group $G$ and consider the conjugation action of $G$
on itself thus inducing an action of $G/Z(G)$ on $G$ given by
$$Inn(G) = G/Z(G) \to Aut(G).$$ Let $\Delta(G)$ denote the associated
extension

\[
\begin{CD}
1 @>{}>> G @>{i}>> \Delta(G) @>{p}>>  G/Z(G) @>{}>> 1
\end{CD}
\] obtained from the conjugation action of $G/Z(G)$ on $G$.

Notice that $SL(n,A)$ acts on the full matrix ring $M(n,A)$ in two
ways recorded next where $M \in M(n,A)$ and $\gamma, y \in SL(n,A)$.

\begin{enumerate}
  \item $(1,y)(M) = yM$ and
  \item $(\gamma,1)(M) = \gamma\cdot(M)\cdot \gamma^{-1}$.
\end{enumerate}

Then define $$(\gamma, y)(M) = \gamma\cdot (y\cdot M) \cdot
\gamma^{-1}.$$

\begin{lem}\label{lem:left.action}
The formula $$(\gamma, y)(M) = \gamma\cdot (y\cdot M) \cdot
\gamma^{-1}$$ for $M \in M(n,A)$ and $\gamma, y \in SL(n,A)$
specifies a left action of $\Delta(SL(n,A))$ on $M(n,A)$.

\end{lem}

\

Assume Lemma \ref{lem:left.action} for the moment.

\begin{thm} \label{thm:actions}

The formula $$\rho((\gamma, y))(M) = \gamma\cdot (y\cdot M) \cdot
\gamma^{-1}$$ for $M \in M(n,A)$ and $\gamma, y \in SL(n,A)$ induces
a faithful representation
$$\rho\colon \Delta(SL(n,A)) \to GL(n^2,A).$$

\end{thm}

\

The theorem has an elementary, immediate consequence.
\begin{cor} \label{cor:G.L.N.squared}

If $G$ is a group with trivial center and is a subgroup of
$SL(n,A)$, then the split extension

\[
\begin{CD}
1 @>{}>> G @>{i}>> \Delta(G) @>{p}>>  G/Z(G) = G @>{}>> 1
\end{CD}
\]

\

\noindent where G acts on itself by conjugation is a subgroup of $GL(n^2,A).$
\end{cor}

\

The proof of Theorem \ref{thm:actions} is given next.
\begin{proof}
First notice that by Lemma  \ref{lem:left.action} the function
$\rho$ is a homomorphism.

If $(\gamma, y)$ is in the kernel of $\rho$ then $$\rho((\gamma,
y))(M) = M$$ for all $M \in M(n,A)$. Let $M=1_{n}$ the
multiplicative identity element in $M(n,A)$. Then
$$\gamma\cdot (y\cdot 1_{n}) \cdot \gamma^{-1}= 1_{n}$$ implies that
$y = 1_{n}.$

Thus assume that $(\gamma, 1_{n})$ is in the kernel of $\rho$. Hence
$$\gamma(M)\gamma^{-1} = M $$ for all $M \in M(n,A)$, and $\gamma$ is in the
center of $PSL(n,A)$ which, by definition is trivial. The Theorem
follows.

\end{proof}

The proof of Lemma \ref{lem:left.action} is given next.
\begin{proof} Let $(\alpha,x)$ and $(\beta,y)$ denote elements
in the semi-direct product $\Delta(SL(n,A))$.

Then the following hold for $M \in M(n,A)$.

\begin{enumerate}
  \item $(\alpha, x)(\beta,y) = (\alpha \beta, \beta^{-1} x \beta y)$
   \item $(\alpha \beta, \beta^{-1} x \beta y)(M) = \alpha \beta(\beta^{-1}x\beta yM)\beta^{-1}\alpha^{-1}
   = \alpha (x\beta yM)\beta^{-1}\alpha^{-1} $
\item $ (\alpha, x)((\beta,y)(M)) =(\alpha, x)(\beta  yM
\beta^{-1})= \alpha( x \beta(yM)\beta^{-1}) \alpha^{-1}$
\end{enumerate} Since the two formulas agree, the Lemma follows.

\end{proof}

\end{document}